\documentclass[12pt,reqno,final]{article}
\usepackage{clz}
\hypersetup{pdfpagemode=UseOutlines,colorlinks=false,pdfpagelayout=SinglePage,pdfstartview=FitH,bookmarksopen=true}
\begin{document}
\title{Lie Rinehart Bialgebras for Crossed Products}
\author{\textsc{Zhuo Chen} \\
{\small LMAM, Peking University \& CEMA, Central University of Finance and Economics,} \\
{\small \href{mailto:chenzhuott@gmail.com}{\texttt{chenzhuott@gmail.com}}}
\and \textsc{Zhang-Ju Liu} \\
{\small  Department of Mathematics and LMAM, Peking University} \\
{\small \href{mailto:liuzj@pku.edu.cn}{\texttt{liuzj@pku.edu.cn}}}
\and \textsc{De-Shou Zhong}
\\
{\small Department of
Economics, CYU}\\
{\small
\href{mailto:zhongdeshou@263.net}{\texttt{zhongdeshou@263.net}}} }
\date{}
\maketitle

\begin{abstract} In this paper, we study Lie Rinehart bialgebras,
the algebraic generalization of Lie bialgebroids. More precisely, we
analyze the structure of Lie Rinehart bialgebras for  crossed
products induced by actions of Lie algebras on $\Kt$.



\end{abstract}


\section*{Introduction}
\addcontentsline{toc}{section}{Introduction}


The concept of Lie Rinehart algebras was introduced in
\cite{MR0154906} as an abstract algebraic treatment of the category
of Lie algebroids \cite{MR1037400,MR1325261}  and were investigated
further in many texts
\cite{MR0197492,MR1058984,MR1625610,MR1696093}.
For the details and history of the notion of Lie Rinehart algebras,
one may refer to an expository paper of Huebschmann
\cite{MR2075590}.  Lie Rinehart algebras may be seen as an algebraic
generalization of the notion of Lie algebroid in which the space of
sections of the vector bundle is replaced by a module over a ring,
vector fields by derivations of the ring and so on. Any attempt to
extend Lie algebroid theory to singular spaces leads to Lie Rinehart
algebras. The reader interested in Lie algebroids and groupoids is
referred to Mackenzie's new book \cite{MR2157566}  (see also
\cite{MR1747916} [Chapters 8 and 12] \cite{MR896907} [I and III])
for background information.

In this paper, we always assume that $\cawu$ is a commutative,
associative algebra over $K$ ($\mathbb{R}$ or $\mathbb{C}$) with a
unit. A Lie Rinehart algebra is an $\cawu$-module that admits a
$K$-Lie bracket and an action on $\cawu$ (called anchor), which are
compatible in a certain sense (see Definition~\ref{Def:DiffLieAlg}).

The notion of Lie Rinehart bialgebras is first introduced by
Huebschmann \cite{MR1764437}. It is derived from the notion of Lie
bialgebras, introduced by Drinfel'd in \cite{MR688240}, and the
notion of Lie bialgebroids introduced by Mackenzie and Xu in
\cite{MR1262213} as the infinitesimal objects associated to a
Poisson groupoid. Lie algebroids are generalized tangent bundles,
while Lie bialgebroids can be considered as generalizations of both
Poisson structures and Lie bialgebras. Roughly speaking, a Lie
bialgebroid is a Lie algebroid $A$ whose dual $A^*$ is also equipped
with a Lie algebroid structure, which is compatible in a certain
sense with that of $A$. This compatibility condition can be
expressed equivalently in terms of the pair $(A,d_*)$, where
$d_*:\Gamma(\wedge^\bullet A)\to\Gamma(\wedge^{\bullet +1} A)$ is
the differential operator inducing the Lie algebroid structure on
$A^*$. Analogously, a Lie Rinehart bialgebra $(E,\sigma)$ is a Lie
Rinehart algebra $E$  endowed with a graded operator
$\sigma:\wedge^k_{\cawu}E\to \wedge^{k+1}_{\cawu}E$ satisfying
$\sigma^2=0$ and a condition similar to Drinfel'd's cocycle
condition (see Definition~\ref{Def:DiffLieBialgebra}). A similar
concept, namely generalized Lie bialgebras, was defined in
\cite{MR1626622}, where the two Lie algebras are not dual in the
usual sense.

In \cite{MR688240}, Drinfel'd classified Lie bialgebras
successfully. The classification of Lie Rinehart bialgebras is, in
our humble opinion, more challenging. However, if we restrict our
attention to the special class of  crossed products  (first
introduced by Malliavin in \cite{MR959753}), we can get lots of
information about their algebraic structure and Lie Rinehart
bialgebra structures over them. The correspondence of crossed
products of Lie algebroids  are known as action Lie algebroids or
transformation Lie algebroids.


The purpose of this paper is twofold. First, we investigate the
structure of actions of a Lie algebra on the ring of polynomials
$\Kt$  and we describe the special features of a $\Kt$-crossed
product\footnote{In this paper we only deal with actions of Lie
algebras on $\Kt$. But in many situations we have to consider the
fraction field $\Ktq$. We expect to get a classification of Lie
algebra actions on $\Ktq$ and further results concerning Lie
Rinehart bialgebras.}. Second, we classify all possible Lie Rinehart
bialgebras $(\Ktq\otimes \LieG,d_*)$ ($\Ktq$ denoting the fraction
field of $\Kt$) in which $\Ktq\otimes \LieG $ is a crossed product
coming from an action of $\LieG$ on $\Kt$. It turns out that our
classification is very similar to the results of \cite{MR1930079}
and \cite{MR2166121}, i.e. the operator $d_*$ (which determines the
dual Lie Rinehart algebra structure) is the sum
$[\Lambda,\cdot]+\Omega$ of a bivector $\Lambda$ and some cocycle
$\Omega$. As in \cite{MR2166121}, we call the data
$(\Lambda,\Omega)$ a compatible pair. In the particular case that
$\LieG$ is semisimple, the Lie Rinehart bialgebra structure is
related to the so-called $\varepsilon$-dynamical $r$-matrices. In
fact, it is a special case of the dynamical $r$-matrices coupled
with Poisson manifolds introduced in \cite{MR1930079}.

The paper is organized as follows.

Section  \ref{Sec:LieRinehartalgebras} reviews Lie Rinehart
algebras, Lie Rinehart bialgebras and   crossed products. Most
importantly, we recall some fundamental properties of the Schouten
bracket and Gerstenhaber algebras.

Section \ref{Sec:actionLieRinehartAlgebraBialgebraKt} builds on the
foundations laid forth in  \cite{MR2283149}, namely the
classification of actions of a finite dimensional Lie algebra on
$\Kt$. Its main result is Theorem~\ref{Thm:qeyqrahg}, which asserts
that any  $\Kt$-crossed product is an extended one.

Section \ref{Sec:actionLieRinehartktBi} is devoted to the
classification of  Lie Rinehart coalgebras for $\Kt$-crossed
products. The main results are
Theorems~\ref{Thm:compatiblepairdstar} and
\ref{Thm:compatiblepairdstar2}, which can be summarized as follows.
If $\Kt\otimes \LieG $ is nontrivial, then the differential operator
of any bialgebra $( \Ktq \otimes \LieG,d_*)$  decomposes as
$d_*=[\Lambda,\cdot]+\Omega$, where $\Lambda$ is a bivector of
$\Ktq\otimes \LieG  $  and $\Omega$ is a map from $\LieG$ to
$L^2=L\wedge_{\Ktq}L$   with $L$ denoting the kernel of
$\theta:~\Ktq\otimes \LieG   \to \Ktq$. The data $(\Lambda,\Omega)$
is called a compatible pair of $\Ktq\otimes \LieG $.

Section \ref{Sec:SemisimpeactionLRBaglebras} details the special
properties enjoyed by the data $(\Lambda,\Omega)$ in the particular
case that $\LieG$ is a semi-simple Lie algebra. Notice that, in this
case, any nontrivial action of $\LieG$ merely comes from
$\mathfrak{sl}(2,K)$, which must be an ideal of $\LieG$. The
conclusion is that the corresponding Lie Rinehart bialgebras
$(\Kt\otimes \LieG ,d_*)$ can be characterized by an
$\varepsilon$-dynamical $r$-matrix $\Lambda$ such that
$d_*=[\Lambda,\cdot]+\varepsilon \mathcal{D}$. Here
$\mathcal{D}:\Kt\otimes \wedge^k\LieG \to \Kt\otimes
\wedge^{k+1}\LieG$ is a fixed operator and $\varepsilon$ is a
constant number in $K$ (see Definition~\ref{Def:epsilonDYBE} and
Theorem~\ref{Thm:LieGOtimesKtBialgebraStrucuture}). In this case,
the data
$(\Lambda+\varepsilon\tau,\varepsilon(\mathcal{D}-[\tau,\cdot]))$ is
a compatible pair, where $\tau\in \Kt \otimes
\wedge^2\mathfrak{sl}(2,K) $ is an $\varepsilon$-dynamical
$r$-matrix with $\varepsilon=-1$ (see
Proposition~\ref{Pro:eteqwt;eral;jkg}).

\paragraph{Acknowledgments}
Besides CPSF(20060400017), the present work was completed while Zhuo
Chen was visiting the mathematics department at Penn State with
support from its Shapiro fund. D.-S. Zhong is grateful to Hongkong
University for its hospitality. The authors thank Johannes
Huebschmann and the anonymous referee for many useful comments.

\section{Lie Rinehart (Bi-)Algebras}\label{Sec:LieRinehartalgebras}

Let $\cawu$ be a unitary commutative algebra over $K=\Real$ or
$\Comp$. A derivation of $\cawu$ is a $K$-liner map, $\delta$:
$\cawu\lon \cawu$, satisfying the Leibnitz  rule
$\delta(ab)=\delta(a)b+a\delta(b)$ for all $  a,b\in\cawu$. The
$\cawu$-module $\Dera$ of all derivations of $\cawu$ is  a $K$-Lie
algebrac under the   commutator $[\delta,\lambda]\defbe
\delta\circ\lambda-\lambda\circ\delta$. In particular, for the ring
of polynomials $\Kt$,  any $\delta\in\mathrm{Der}\Kt$ has the form
$\delta=f\frac{d}{dt}$, where $f$ is uniquely determined by
$f=\delta(t)$. Under this expression, we have
$$[f\frac{d}{dt},g\frac{d}{dt}]=(fg'-f'g)\frac{d}{dt},$$
where $f'=\frac{d}{dt}f$. Thus, with the  Lie bracket:
$[f,g]=fg'-f'g$ for all $f,g\in \Kt$, we can identify $\Kt$ with
$\mathrm{Der}\Kt$.


\begin{defn}\label{Def:DiffLieAlg}
A Lie Rinehart algebra is a pair $(\cawu,\ME)$ where $\ME$ is both
an $\cawu$-module as and a $K$-Lie algebra such that
\begin{itemize}
\item[1)]
 there is a Lie algebra morphism $\theta$:
$\ME\lon  \Dera$ {\rm (called the anchor of $\ME$)} which is also a
morphism of $\cawu$-modules;
\item[2)]$ [ X_1,a X_2] = a[ X_1, X_2]+\theta(
X_1)(a) X_2$,  $\forall  X_i\in \ME, a\in \cawu$.
\end{itemize}
\end{defn}
In  case that $\cawu$ is fixed  we just say that $\ME$ is a Lie
Rinehart algebra.

 \vskip 0.2cm \noindent\textbf{$\bullet$  Crossed Products}

Let $(\LieG, [\cdot,\cdot])$ be a   $K$-Lie algebra. We denote the
$\cawu$-module $\cawu\otimes_K\LieG $ by $\LieGOtimesa$ and write an
element $a\otimes X $ as $aX$. An action of $\LieG$ on $\cawu$ means
a Lie algebra morphism $\theta:\LieG\lon \Dera$. By  the same symbol
 $\theta: \LieGOtimesa\lon \Dera$ to denote the $\cawu$-module
morphism extended from  this action, then  we have an induced
bracket defined on $\LieGOtimesa$:
$$
\ActionBracket{aX}{bY}\defbe ab[X,Y]+a(\theta(X)b)Y-b(\theta(Y)a)X,
\quad \forall a,b\in\cawu, X,Y\in \LieG
$$
such that  $\LieGOtimesa $ is a Lie Rinehart algebra.
\begin{defn}[\cite{MR959753}]The triple $(\LieGOtimesa,\ActionBracket{\cdot}{\cdot},\theta)$
is called a crossed product, generated by $\LieG$ via the action
$\theta$, which is said to be nontrivial (resp. trivial) if $\theta$
is nontrivial (resp. trivial).
\end{defn}
 The
reader is recommended to compare with the so-called ``action Lie
algebroids'' or ``transformation algebroids''
\cite{MR896907,MR1409974,MR1747916} to understand the geometric
background of crossed products.

To introduce a special kind of   crossed products, we define firstly
the notion of a derivation of an $\cawu$-Lie algebra.
\begin{defn}\label{Def:derivationofL}Let $(L,[\cdot,\cdot]_L)$ be  an $\cawu$-Lie
algebra, i.e., a Lie Rinehart algebra  with a trivial anchor. A
derivation of $L$ is a pair $(\huaD,\delta)$, where $\huaD: L\lon L$
is a $K$-linear operator, $\delta\in\Dera$ and they satisfy the
conditions:
\begin{eqnarray*}
\huaD[l_1,l_2]_L&=&[\huaD l_1,l_2]_L+[l_1,\huaD l_2]_L,\quad\forall
l_1,l_2\in L,\\
\huaD(a l)&=&\delta(a)l + a \huaD l,\quad\forall a\in \cawu, l\in L.
\end{eqnarray*}
\end{defn}

\begin{prop}
For a derivation $(\huaD,\delta)$
  as above, $\cawu\oplus L$ has a natural  crossed product  structure  given by:
\begin{eqnarray*}
\theta ( a, l) &=& a\delta,\\
{[(a_1,l_1),(a_2,l_2)]} &=& ({a_1} \delta {(a_2)}- a_2 \delta
{(a_1)}, [l_1,l_2]_L+ a_1 \huaD l_2- a_2 \huaD l_1)\,,
\end{eqnarray*}
for all $(a,l)$, $(a_i,l_i)\in \cawu\oplus L$.
\end{prop}
We will call the Lie Rinehart algebra constructed in this way an
\emph{extended   crossed product} of $L$ via a derivation
$(\huaD,\delta)$, and it will be denoted by $\cawu
\ltimes_{(\huaD,\delta)} L$. It is said to be \emph{nontrivial}, if
$\delta$ is not zero.

\begin{defn}Let $(\huaD_i,\delta_i)$ ( $i=1,2$ ) be two derivations of
 $\cawu$-Lie algebras $(L_i,[\cdot,\cdot]_i)$ respectively. They are said to be
equivalent, written $(\huaD_1,\delta_1;L_1)\sim
(\huaD_2,\delta_2;L_2)$ if there exists an $\cawu$-Lie algebra
isomorphism $\Phi: L_1\lon L_2$, an invertible element $a_0\in\cawu$
and $l_0\in L_2$ such that
\begin{eqnarray}\nonumber
\left\{
\begin{array}{r@{\quad =\quad}l}
\huaD_1 & a_0 \Phi\inverse\huaD_2\Phi +
[\Phi\inverse(l_0),\cdot]_1\,,\\
\delta_1  & a_0\delta_2\,.
\end{array}
\right.
\end{eqnarray}
\end{defn}
It is obvious that ``$\sim$'' is equivalence relation and the
following proposition  is  easy to be verified.
\begin{prop}\label{Pro:extendedLRAClassification}
Assume that $\cawu$ has no zero-divisors. Let $(\huaD_i,\delta_i)$ (
$i=1,2$ ) be respectively derivations of the $\cawu$-Lie algebras
$(L_i,[\cdot,\cdot]_i)$  and assume that $\delta_i\neq 0$. Then
$\cawu \ltimes_{(\huaD_i,\delta_i)} L_i$ ($i=1,2$) are isomorphic if
and only if $(\huaD_1,\delta_1;L_1)\sim (\huaD_2,\delta_2;L_2)$.
\end{prop}

\vskip 0.2cm \noindent\textbf{$\bullet$ The Schouten bracket and
Gerstenhaber algebras}

 A Gerstenhaber algebra consists of a triple
$(A=\sum_{i\in\mathbb{Z}}A^i,\wedge,[\cdot,\cdot])$ such that
$(A,\wedge)$ is a graded commutative associative algebra over $K$
and $(A=\sum_{i\in\mathbb{Z}}A^{(i)},[\cdot,\cdot])$ is a graded Lie
algebra, where $A^{(i)}=A^{i+1}$, such that $[a,\cdot]$ is a
derivation with respect to $\wedge$ of degree $(i-1)$ for any $a\in
A^{i}$.

 It is
shown in \cite{MR1077465} that  a  Lie Rinehart algebra $E$
corresponds a Schouten algebra $\wedgea^\bullet E$, which is, in
fact, a Gerstenhaber algebra \cite{MR1361447} (see also Theorem 5 in
\cite{MR1187279}). The Schouten bracket is a $K$-bilinear bracket
$[\cdot,\cdot]$: $\wedgea^k E\times \wedgea^l E\lon$
$\wedgea^{k+l-1} E$ such that $(\wedgea^\bullet
E,\wedgea,[\cdot,\cdot])$ forms a Gerstenhaber algebra such that:
\begin{itemize}
\item[a.] It coincides with the original Lie bracket on $E$.
\item[b.] $[x,f]=\theta(x)f$, $\forall f\in \cawu$, $x\in E$.
\item[c.] It is a derivation in the graded sense, i.e.,
\begin{equation}\nonumber
[x,y\wedgea z]=[x,y]\wedgea z+(-1)^{(\abs{x}-1)\abs{y}}y\wedgea
[x,z],
\end{equation}
\end{itemize}
where $x\in \wedgea^{\abs{x}} E$, $y\in \wedgea^{\abs{y}} E$, $z\in
\wedgea^\bullet E$.

Conversely, the axioms of a Gerstenhaber algebra
$(A=\sum_{i\in\mathbb{Z}}A^i,\wedge,[\cdot,\cdot])$ naturally imply
that $(A^0,A^1)$ is a Lie Rinehart algebra, such that $\theta(x)f=
[x,f]$, for each $x\in A^1$ and $f\in A^0$.

Let $ E $ be a Lie Rinehart algebra  and $F$ an $\cawu$-module. By
saying a representation of $E$ on $F$, we mean an $\cawu$-map:
$E\times F\lon F$, $x\times s\mapsto x.s$, satisfying the following
axioms:
\begin{eqnarray*}
(fx).s &=& f(x.s);\\
x.(fs) &=& f(x.s)+\theta(x)(f)s;\\
x.(y.s)-y.(x.s) &=& [x,y].s,\ \ \forall s\in F, x,y\in E, f\in
\cawu.
\end{eqnarray*}
An $\cawu$-map $\Omega: E\lon F$ is called a $1$-cocycle, if
\begin{equation}\label{Eqt:Cocycle}
\Omega[x,y]=x.\Omega(y)-y.\Omega(x),\ \ \forall x,y\in E.
\end{equation}

For example, let $L=\Ker\theta$, which is clearly an $\cawu$-module,
as well as an ideal. We define the \emph{adjoint representation} of
$E$ on $L$: (or on $\wedgea^k L$ in the sense of the Schouten
bracket, for some $k\geq 2$ )
$$x.l\defbe [x,l],\ \ \forall x\in E,\ l\in L \ (\mbox{or }\wedgea^k L).$$

\vskip 0.2cm\noindent\textbf{$\bullet$ Lie Rinehart bialgebras}

A differential Gerstenhaber algebra is a Gerstenhaber algebra
equipped with a derivative operator $\sigma$, called the
differential, which is of degree $1$ and square zero. It is called a
strong differential Gerstenhaber algebra if $\sigma$ is also a
derivation of the graded Lie bracket \cite{MR1675117}. We recall a
similar concept, namely the Lie Rinehart bialgebras, first
introduced by Huebschmann \cite{MR1764437}.

A graded operator (of degree $1$) on $\wedgea^\bullet E $ is a
$K$-linear operator $\sigma$: $\wedgea^k E \lon \wedgea ^{k+1} E$
satisfying
\begin{equation}\nonumber
\sigma(x\wedgea y)=\sigma x\wedgea y+(-1)^{\abs{x}}x\wedgea \sigma
y, \ \ \forall x\in \wedgea^{\abs{x}} E,y\in \wedgea^\bullet E.
\end{equation}
Let $\sigma$ be a graded operator, then it induces two structures on
$\EXinga$ = $\Hom_{\cawu}(E,\cawu)$, the $\sigma$-anchor
$\theta_\sigma$ and $\sigma$-bracket $[\cdot,\cdot]_\sigma$, such
that
\begin{eqnarray*}
\theta_\sigma(\xi)f&=&<\sigma f,\xi>\\
<[\xi,\eta]_\sigma,x> &=&- <\sigma
x,\xi\wedgea\eta>+\theta_\sigma(\xi)<x,\eta>
-\theta_\sigma(\eta)<x,\xi>.\ \
\end{eqnarray*}
The following proposition is a fundamental criterion.
\begin{prop}[\cite{MR1077465}]\label{Pro:PreLieAlgExteriorOperatorStar}
Equipped with the two structures given by the graded operator
$\sigma$, $\EXinga$ is a Lie Rinehart algebra if and
only if $\sigma^2=0$. 
\end{prop}

\begin{defn} \label{Def:DiffLieBialgebra}
Let $\sigma$: $\wedgea^k E \lon \wedgea ^{k+1} E$ be a graded
operator of degree $1$. If $\sigma^2=0$ and
\begin{equation}\label{Eqt:LieBialgebradXingxf}
\sigma [x,y]=[\sigma x, y]+(-1)^{(\abs{x}+1)}[x,\sigma y],\quad
\forall x\in \wedgea^{\abs{x}} E,y\in \wedgea^\bullet E.,
\end{equation}
$(E, \sigma)$ is called a  Lie Rinehart bialgebra.
\end{defn}
Thus for a strong differential Gerstenhaber algebra
$(\wedgea^\bullet E,\wedgea,[\cdot,\cdot])$ with differential
$\sigma$, $(E,\sigma)$ is naturally a Lie Rinehart bialgebra.   We
omit the proofs of the following three propositions since they are
straightforward.

\begin{prop}\label{Lem:sigmaextend} Let $\sigma$ be
a graded operator of degree $1$ which is also square zero. Suppose
that
\begin{equation}\label{Eqt:LieBialgebradXingxy}
\sigma[x,y]=[\sigma x,y]+[x ,\sigma y], \ \ \ \ \forall x,y \in E.
\end{equation}
\begin{itemize}
\item[1)] If $E$ is \emph{nondegenerate}, i.e., for $x\in E$,
$$
x\wedgea y=0,\ \forall y\in E \quad\mbox{ implies }\quad x=0,
$$
then
\begin{equation}\nonumber
\sigma [x,f]=[\sigma x,f]+[x,\sigma f],\ \ \forall x\in E,\
f\in\cawu.\end{equation}

\item[2)] If $E$ is nondegenerate and
\emph{faithful}, i.e., for $a\in\cawu$,
$$
ax=0,~ \forall x\in E \quad\mbox{ implies }\quad a=0,
$$
then  (\ref{Eqt:LieBialgebradXingxf}) holds.

\item[3)] Define a bracket
\begin{equation}\nonumber
\PoissonBracket{a}{b}=[\sigma a,b], \
\forall a,b\in \cawu.
\end{equation} If $E$
is nondegenerate, then the algebra
$(\cawu,\PoissonBracket{\cdot}{\cdot})$ is a Leibnitz algebra.
Moreover, the bracket is  skew-symmetric if $E$ is faithful, then
$\cawu$ is  a Poisson algebra in this case.

\item[4)]If $E$ is nondegenerate and faithful, then $(\wedgea^\bullet E$,
$\wedgea$, $[\cdot,\cdot])$ is a strong differential Gerstenhaber
algebra
with the differential $\sigma$. 
\end{itemize}
\end{prop}

We are going to discover further properties of Lie Rinehart
bialgebras with some additional conditions. First we review a
special kind of Lie Rinehart bialgebra which generalized the method
that a Poisson tensor $\pi$ on a manifold gives a $\pi$-bracket for
the $1$-forms. This method is also referred as the dualization of a
Lie Rinehart algebra (see Kosmann and Magri's definition in
\cite{MR1077465}).


If $\Lambda$ is a bivector of $E$, i.e., $\Lambda\in \wedgea^2 E$,
then clearly the operator $\sigma=[\Lambda,\cdot]$ satisfies
condition (\ref{Eqt:LieBialgebradXingxf}). When $\sigma^2=0$, or
equivalently, $[[\Lambda,\Lambda],\cdot]=0$, then $(E,\sigma)$ is a
Lie Rinehart bialgebra. For objects of this type, we call them
\emph{coboundary} (or exact) ones \cite{MR1371234}. Especially, when
$[\Lambda,\Lambda]=0$, we call $\Lambda$ a Poisson bivector and
$(E,[\Lambda,\cdot])$ a \emph{triangular} Lie Rinehart bialgebra.
Several examples will be given later after Theorem
\ref{Thm:trivialactionbialgebra}.


Let $E$ be a Lie Rinehart algebra and $\Lambda$ a bivector. We will
use the symbol $\Lambda^\sharp$ to denote the contraction map
$\EXinga\lon E$, defined by  $\Lambda^\sharp(\phi)=\phi\lrcorner
\Lambda $ (this is legal since $\Lambda$ can be expressed as a
finite sum $\sum a_i \wedge b_i $, and thus $\phi\lrcorner
\Lambda=\phi(a_i)b_i-\phi(b_i)a_i$). The operation
\begin{eqnarray*}
[\phi,\psi]_\Lambda &=& d<\Lambda,\phi\wedge\psi> +
{\Lambdasharp(\phi)}\lrcorner d\psi- {\Lambdasharp(\psi)}\lrcorner
d\phi,
\end{eqnarray*}
for $\phi$, $\psi\in \EXinga$, is called the
\emph{$\Lambda$-bracket} on $\EXinga$. Equivalent expressions are
given as follows.
\begin{eqnarray}
\nonumber
&&<[\phi,\psi]_\Lambda,x>\\\label{Eqn:phipsiLambdabracket2}
&=&<[x,\Lambda],\phi\wedge\psi>+\theta({\Lambdasharp(\phi)})<\psi,x>-
\theta({\Lambdasharp(\psi)})<\phi,x>,
\end{eqnarray}
for any two $\phi,\psi\in \EXinga$. We omit the proof of these
relations. By (\ref{Eqn:phipsiLambdabracket2}), one is able to get
the following proposition.
\begin{prop}
Let $E$ be a Lie Rinehart algebra and $\Lambda$  a bivector of $E$
such that $[[\Lambda,\Lambda],\cdot]=0$. Then, for the Lie Rinehart
bialgebra $(E,[\Lambda,\cdot])$, the corresponding Lie Rinehart
algebra $\EXinga$ given by Proposition
\ref{Pro:PreLieAlgExteriorOperatorStar}, has the anchor map
$\theta_*=\theta\circ\Lambdasharp$ and bracket
$[\cdot,\cdot]_*=[\cdot,\cdot]_\Lambda$. 
\end{prop}


\section{Classification of  $\Kt$-Crossed Products}
\label{Sec:actionLieRinehartAlgebraBialgebraKt}

The algebra we study in this section is always assumed to be the
ring of polynomials $\cawu=\Kt$ and $\LieG$ is assumed to be a
finite dimensional Lie algebra over $K$. First recall that any Lie
algebra $\LieG$ admits a unique maximal solvable ideal of $\LieG$,
denoted by $\RadLieG$, and called the \textbf{radical}, or Jacobson
root. The famous Levi decomposition of a Lie algebra is expressed as
$\LieG=\RadLieG\rtimes \frkm$, where $\frkm$ is a semisimple Lie
subalgebra {\rm (known as the \textbf{Levi subalgebra} of $\LieG$,
which is not necessarily unique \cite{MR746308})}.

Next we quote the following result   in \cite{MR2283149} as the
first step of the classification of  crossed products
$\LieGOtimesaKt$, where all actions of an arbitrary Lie algebra
$\LieG$ on $\Kt$ are classified into three types according to
$\Rank(\theta)$, i.e., the dimension of $\Img (\theta)$.

\begin{thm}\cite{MR2283149}\label{Thm:Main}
Let $\LieG$ be a Lie algebra. Let $\RadLieG$ be its radical and
 $\frkm$  a Levi subalgebra. If $\theta: \LieG\lon \Kt$ is a
 nontrivial action,
then $\Rank(\theta)\leq 3$ and the action has following three
possible types:
\begin{itemize}
\item\quad{\rm Type 1:} \quad$\Rank(\theta)=1$. In this case,
$\theta|_\frkm=0$, and there exists a polynomial $h\in\Kt$ and a
linear function $\lambda\in\LieGXing$ (both nonzero), such that
\begin{equation}\label{Eqt:thetatype1}
\theta(X)=\lambda(X)h,\quad \forall X\in\LieG.
\end{equation}
\item\quad{\rm Type 2:} \quad$\Rank(\theta)=2$. In this case,
$\theta|_\frkm=0$, and there exist a nonnegative integer $m\neq 1$,
a constant $b\in K$, and two linearly independent
$\lambda,\mu\in\LieGXing$, such that
\begin{equation}\label{Eqt:thetatype2}
\theta(X)=\lambda(X)(t+b)^m+\mu(X)(t+b),\quad \forall X\in\LieG.
\end{equation}

\item\quad{\rm Type 3:} \quad$\Rank(\theta)=3$. In this case,
one is able to decompose $\frkm=\frks\oplus \frkm_0$, where
$\frks\cong \sltwoK$, $\frkm_0$ is a semisimple Lie subalgebra such
that $\Ker (\theta)=\RadLieG\rtimes\frkm_0$. Moreover, one is able
to find a standard base $X_0$, $X_1$, $X_2\in \frks$, such that
$\theta(X_0)=1$, $\theta(X_1)=t$, $\theta(X_2)=t^2$.
\end{itemize}
\end{thm}


As a standard base of $\sltwoK\subset \gl(2,K)$, three vectors
\begin{equation}\label{Eqt:HEPEMofsl2}
\HH={\half} {\begin{pmatrix}
  {1} & {0} \\
  {0} & {-1}
\end{pmatrix}},
\ \ \EEP=
\begin{pmatrix}
  0 & -1 \\
  0 & 0
\end{pmatrix},\ \
\EEM=
\begin{pmatrix}
  0 & 0 \\
  1 & 0
\end{pmatrix}.
\end{equation}
are related by $ [\HH,\EEP ]=\EEP $, $[\HH,\EEM ]=-\EEM $, $[\EEP
,\EEM ]=-2\HH$.

\begin{ex}
\label{Ex:sl2KstandardMorphism} We define a special $\theta_s$
(called the \emph{standard action}) of $\sltwoK$ on $\Kt$:
$$
\left\{
\begin{array}{r@{\quad=\quad}l}
\theta_s(\HH)& t,\\
\theta_s(\EEP )& t^2,\\
\theta_s(\EEM )& 1.
\end{array}
\right.
$$
It is typically an action of Type 3.
 \end{ex}

\begin{prop}\label{Prop:type123details} With the same assumptions as
in Theorem \ref{Thm:Main}, then
\begin{itemize}
\item[1)] $\LieG$ admits an action of Type 1 if and only if
$[\LieG,\LieG]\subsetneq \LieG$. In this case,
(\ref{Eqt:thetatype1}) defines an action if and only if
$[\LieG,\LieG]\subset \Ker\lambda$.
\item[2)] $\LieG$ admits an action of Type 2 if and only if there
exist two independent vectors $x_0$, $y_0\in\LieG$ and an ideal
$S\subset \LieG$ such that
\begin{itemize}
\item[2.1)] $\LieG= S\oplus \pairing{x_0}\oplus \pairing{y_0}$;
\item[2.2)] $[\LieG,\LieG]\subset S\oplus \pairing{x_0}$;
\item[2.3)] $[x_0,y_0]+(m-1)x_0\in S$, for some nonnegative integer $m\neq
1$.
\end{itemize}
In this case, by setting
$$\lambda|_{S\oplus\pairing{y_0}}=0,\quad\mbox{and}\quad \lambda(x_0)=1,$$
$$\mu|_{S\oplus\pairing{x_0}}=0,\quad\mbox{and}\quad \mu(y_0)=1,$$
 Equation (\ref{Eqt:thetatype2}) defines an action.
\item[3)] $\LieG$ admits an action of Type 3 if and only if
$\LieG$ is not solvable and the Levi subalgebra of $\LieG$ admits
$\sltwoK$ as an ideal.
\end{itemize}
\end{prop}
\begin{proof} The statements (1) and (3) are direct consequences of Theorem
\ref{Thm:Main}. We elaborate on (2). In fact, a simple calculation
shows that, the map $\theta$ defined by (\ref{Eqt:thetatype2}) is an
action if and only if for $X,Y\in \LieG$,
\begin{eqnarray}\nonumber
\mu([X,Y])&=&0,\\\nonumber 
\lambda([X,Y])&=&(1-m)(\lambda(X)\mu(Y)-\lambda(Y)\mu(X)).
\end{eqnarray}
In this case,  $S=\Ker\mu\cap \Ker\lambda$ is an ideal. We select
$x_0,y_0\in \LieG-S$, satisfying $\lambda(x_0)=\mu(y_0)=1$,
$\lambda(y_0)=\mu(x_0)=0$. Then the three conditions   2.1) $\sim$
2.3) are satisfied. The converse is also obvious.
\end{proof}

Recall Proposition \ref{Pro:extendedLRAClassification} which gives a
classification of nontrivial extended crossed products using the
data $(\huaD,\delta)$. The following theorem claims that any
nontrivial crossed product  $\LieGOtimesaKt$ can be realized as an
extended one. Thus we obtain a classification of such  crossed
products.

\begin{thm}\label{Thm:qeyqrahg}
 For any nontrivial action $\theta:
\LieG\lon \Kt$, the corresponding  crossed product $\LieGOtimesaKt$
is isomorphic to an extended  crossed product $
\Kt\ltimes_{(\huaD,\delta)} L$. The data $L$, $\huaD$ and $\delta$
are respectively specified as follows:
\begin{itemize}
\item[1)]if $\theta$ is of Type 1 defined by (\ref{Eqt:thetatype1}),
then $L=\Kt\otimes \Ker\lambda $, $\huaD=[x_0,\cdot]$.
$\delta=h\frac{d}{dt}$, where $x_0\in\LieG$ satisfies
$\lambda(x_0)=1$;
\item[2)]if $\theta$ is of Type 2 defined by (\ref{Eqt:thetatype2}),
then $L=\Kt\otimes (S\oplus \pairing{x_0-(t+b)^{m-1}y_0}) $,
$\huaD=[y_0,\cdot]$, $\delta=(t+b)\frac{d}{dt}$, where $x_0$, $y_0$
and $S$ are specified by (2) of Proposition
\ref{Prop:type123details};
\item[3)]if $\theta$ is of Type 3,
then $$L=\Kt \otimes (\RadLieG\oplus \frkm_0\oplus
\pairing{X_2-tX_1, X_1-tX_0}), $$ $\huaD=[X_0,\cdot]$ and
$\delta=\frac{d}{dt}$, where $X_0$, $X_1$, $X_2$ is a basis of
$\frks\cong \sltwoK$ declared in (3) of Theorem \ref{Thm:Main}.
\end{itemize}
\end{thm}
\begin{proof} 1) By (1) of Proposition \ref{Prop:type123details},
$[\LieG,\LieG]\subset \Ker\lambda$, $\huaD(L)\subset L$, and hence
we have the conclusion.

2) One only need to check that $\huaD(L)\subset L$.

3) By Theorem \ref{Thm:Main}, we conclude that the Levi subalgebra
$\frkm=\frks\oplus \frkm_0$ and $\theta$ must be trivial on
$\RadLieG\rtimes \frkm_0$. Moreover, $ [X_1,X_2]=X_2$,
$[X_1,X_0]=-X_0$, $[X_2,X_0]=-2X_1$. Thus, the  crossed product
$\Kt\otimes \LieG $ is spanned (over $\Kt$) by:
$$X_0,\quad
A=X_2-tX_1,\quad B=X_1-tX_0\quad\mbox{ and elements in }
\RadLieG\rtimes \frkm_0\,.$$ Clearly, $L=\Kt \otimes (\RadLieG\oplus
\frkm_0\oplus \pairing{A,B}) $ is the kernel of $\theta$:
$\Kt\otimes \LieG \lon \Kt$. And $(\huaD=[X_0,\cdot],\frac{d}{dt})$
is a derivation of $L$. In this way, the Lie Rinehart algebra
$\Kt\otimes \LieG \cong \Kt\ltimes_{(\huaD,\frac{d}{dt})} L$ by
identifying $X_0=1_{{\Kt}}$. \end{proof}

\vskip 0.2cm

Consider $\Ktq$, the fractional filed of $\Kt$ and treat $\Ktq$ as a
$\Kt$-module, as well as for $\LieGOtimesaKtq$. We also have
$\mathrm{~Der}\Ktq\cong \Ktq$. Note that $\LieGOtimesaKtq$ is a
$\Ktq$-crossed product  in an obvious sense.

\begin{prop}\label{Pro:Existencegamma}
 Let $\theta$
be a nontrivial action of $\LieG$ on $\Kt$. Then,
\begin{itemize}
\item[1)] one can find a $\Kt$-module map $\gamma: \Kt\lon
\LieGOtimesaKtq$ such that $ \theta\circ \gamma=Id_{\Kt} $;
\item[2)] the corresponding  crossed product is of extended type, i.e.,
for $L$ being the kernel of $\theta: \LieGOtimesaKtq\lon \Ktq$, and
  $(\huaD,\frac{d}{dt})$   a derivation of $L$, one has
$\LieGOtimesaKtq\cong \Ktq \ltimes_{(\huaD,\frac{d}{dt})} L$.
\item[3)] if $\theta$ is of Type 3, $\gamma$ takes values in
$\LieGOtimesaKt$.
\end{itemize}
\end{prop}
\begin{proof}   We directly construct $\gamma$.
If $\theta$ is of Type 1, then we can find an $X\in
\LieG-[\LieG,\LieG]$ such that $\lambda(X)=1$. In this case, we set
$\gamma(1)= \frac{1}{h} X$. If $\theta$ is of Type 2, then one can
also find an $X\in \LieG-[\LieG,\LieG]$ such that $\mu(X)=1$. In
this case, we set
$$\gamma(1)= \frac{1}{\lambda(X)(t+b)^m+(t+b)} X,$$ whence the first
statement. If $\theta$ is of Type 3, we set $\gamma(1)=X_0$, which
was declared in (3) of Theorem \ref{Thm:Main}. The second statement
is a direct consequence of the first one.
\end{proof}

\vskip 0.2cm

We finally notice the interesting fact that by Proposition
\ref{Prop:type123details}, for any nontrivial action $\theta$ of a
semisimple $\LieG$ on $\Kt$, the effective part of this action
merely comes from $\sltwoK$ and $\theta$ must be of Type 3.
\begin{thm}\label{Thm:sl2KOtimesKtAllEqual}
Let $E=\Kt\otimes \sltwoK $ be a crossed product, with the structure
coming from a nontrivial action $\theta: \sltwoK\lon \Kt$. What ever
$\theta$ is chosen, all these $E$ are isomorphic to each other.
\end{thm}
\begin{proof} It is not hard to see the following fact: for any two nontrivial morphisms
$\theta_1,\theta_2$: $\sltwoK\lon \Kt$, there exists an automorphism
$\Pi$ of $\sltwoK $   such that $\theta_1=\theta_2\circ\Pi$. We
write $E^1$, $E^2$ to indicate the two crossed products. Now we
define an isomorphism $\overline{\Pi}$ from $E^1$ to $E^2$, which
maps $fX$ to $f\Pi(X)$ ($f\in \Kt$, $X\in \sltwoK$). The second
statement is a direct consequence of the first one.
\end{proof}
In general, for a semisimple Lie algebra $\LieG$ which admits an
action of Type 3, it has a  simple ideal isomorphic to $\sltwoK$, so
that $\LieG\cong \sltwoK\oplus \frkm_0$, for some semisimple ideal
$\frkm_0$. Any nontrivial  crossed product $\LieGOtimesaKt$ is
isomorphic to $\Kt\otimes \sltwoK \oplus \Kt\otimes \frkm_0 $, where
$\sltwoK$ has the standard action and $\frkm_0$ has the trivial
action on $\Kt$.

\section{Lie Rinehart Coalgebras  for Crossed Products}
\label{Sec:actionLieRinehartktBi}

It seems that for a non-semisimple Lie algebra $\LieG$, the
structure of a Lie Rinehart bialgebra of crossed product
$(\LieGOtimesaKt,d_*)$ is quite complicated.  We shall discuss the
situation that  $\LieG$ is semisimple in the next section. However,
we can still say something about the operator $d_*$.

 In \cite{MR2166121}, we proved that for a transitive Lie
algebroid $(\huaA,[\cdot,\cdot]_\huaA,\rho)$, the structure of any
Lie bialgebroid $(\huaA,d_*)$ can be characterized by a bisection
$\Lambda\in\Gamma(\wedge^2 \huaA)$ and a Lie algebroid $1$-cocycle,
 $\Omega: \huaA\lon \wedge^2 L$,  with  respect to the adjoint representation of   $\huaA$  on $\wedge^2
 L$,
 where $L=\Ker \rho$ is the
  isotropic  bundle of $\huaA$. Moreover, such a pair is unique up to a gauge term in $\Gamma(\wedge^2 L)$ and the
   differential $d_*$ is decomposed into
    $$
    d_*=[\Lambda,\cdot]_\huaA+\Omega.$$
We will show some similar results of Lie Rinehart bialgebras for
$\LieGOtimesaKt$. 

\begin{defn}\label{Def:1cocycle}Let $\LieG$ be a Lie algebra and let
$\theta$ be an action of $\LieG$ on $\Kt$. For the  crossed product
$(\LieGOtimesaKtq,\ActionBracket{\cdot}{\cdot}, \theta)$, let $L$ be
the kernel of $\theta:~\LieGOtimesaKtq\lon \Ktq$ and
$L^2=L\wedge_{\Ktq} L$. A $K$-linear map $\Omega:~~\LieG\lon L^2$ is
called a $1$-cocycle  if
$$
\Omega[X,Y] =[\Omega(X),Y]+[X,\Omega(Y)],\quad\forall X,Y\in \LieG.
$$
\end{defn}

Such a $1$-cocycle $\Omega$ can be extended as a derivation of the
graded module $\Omega$: $\Ktq\otimes \wedge^{k} \LieG  \lon
\Ktq\otimes \wedge^{k+1} \LieG  $, $k\geq 0$. For $k=0$, it is zero.
For $k=1$, it is simply defined by $fX\mapsto f\Omega(X)$, $\forall
f\in\Ktq$, $X\in \LieG$. For $k> 1$, it is defined by
\begin{eqnarray*}
&&\Omega(u_1\wedge_{\Ktq} \cdots\wedge_{\Ktq}
u_k)\\
&=&\sum_{i=1}^{k}(-1)^{i+1}u_1\wedge_{\Ktq} \cdots\wedge_{\Ktq}
\Omega (u_i)\wedge_{\Ktq} \cdots\wedge_{\Ktq} u_k,\\
&&\ \ \forall u_1,\cdots,u_k\in \LieGOtimesaKtq.
\end{eqnarray*}
One has the following formula:
$$\Omega[u,v] =
[\Omega(u),v] +(-1)^{k+1}[u,\Omega(v)] ,\quad\forall u\in
\Ktq\otimes \wedge^k\LieG,\ \ v\in \Ktq\otimes \wedge^l\LieG.$$

\begin{defn}
With the  assumptions of Definition \ref{Def:1cocycle}, given
$\Lambda \in \Ktq\otimes \wedge^2\LieG $ and a $K$-linear map
$\Omega$: $\LieG\lon {L^2}$, the  pair $(\Lambda,\Omega)$ is called
compatible if $\Omega$ is a 1-cocycle and satisfies
\begin{equation}\nonumber
[\frac{1}{2}[\Lambda,\Lambda]+\Omega(\Lambda) ,\cdot]+\Omega^2=0,
\quad\mbox{as a map } \Ktq\otimes\wedge^2\LieG  \lon
\Ktq\otimes\wedge^3\LieG .
\end{equation}
\end{defn}

If $(\Lambda,\Omega)$ is compatible, then so is the pair
$(\Lambda+\nu,\Omega-[\nu,\cdot])$, for any $\nu\in {L^2}$. Thus,
two compatible pairs $(\Lambda,\Omega)$ and $(\Lambda',\Omega')$ are
called equivalent, written $(\Lambda,\Omega)\sim(\Lambda',\Omega')$,
if $\exists\nu\in{L^2}$, such that $\Lambda'=\Lambda+\nu$ and
$\Omega'=\Omega-[\nu,\cdot]$.

\begin{thm}\label{Thm:compatiblepairdstar}
Let $\LieG$ be a  Lie algebra and let $\theta$ be a nontrivial
action of $\LieG$ on $\Kt$. Then, there is a one-to-one
correspondence between Lie Rinehart bialgebras
$(\LieGOtimesaKtq,d_*)$   and equivalence classes of compatible
pairs $(\Lambda,\Omega)$ such that
\begin{equation}\label{EqtRelationDStarAndPair}
d_*=[\Lambda,\cdot]+\Omega.
\end{equation}
\end{thm}
We first prove the following lemma.
\begin{lem}With the same assumptions, there exists some $\Lambda  \in \Ktq\otimes\wedge^2\LieG $ such that
\begin{equation}\label{temp:5325656}
d_* f=[\Lambda,f],\quad\forall f\in\Ktq.
\end{equation}
\end{lem}
\begin{proof} According to  Proposition
\ref{Pro:Existencegamma} , we can find an element $\gamma(1)\in
\LieGOtimesaKtq$ such that $ \theta\circ \gamma(1)=1$. Then we
define $\Lambda  \in \Ktq\otimes\wedge^2\LieG $ by setting
$$
\Lambda\defbe d_*t \wedge_{\Ktq} \gamma(1).
$$
By (3) of Proposition \ref{Lem:sigmaextend}, there is an
antisymmetric pairing satisfying
$$
\PoissonBracket{f}{g}=[d_*f,g]=f'g'[d_*t,t],\quad\forall f,g\in\Ktq.
$$
Hence it must be zero, i.e., $[d_*f,g]=0$. So we have
\begin{eqnarray*}
&&[\Lambda,f]=-[d_*t,f]\gamma(1)+[\gamma(1),f]d_*t\\
&=&f'\theta\circ\gamma(1)d_*t=f'd_*t=d_*f,\quad\quad\forall
f\in\Ktq,
\end{eqnarray*}
whence the result.\end{proof}

\begin{proof}[Proof of Theorem \ref{Thm:compatiblepairdstar}] Suppose that a Lie Rinehart bialgebra $(\LieGOtimesaKtq,d_*)$ is
given. Then with $\Lambda$ given as in the above lemma, we define
$$
\Omega=d_*-[\Lambda,\cdot].
$$
Equation (\ref{temp:5325656}) implies that  $\Omega$ satisfies
$$\Omega(fu)=f\Omega(u),\quad\;\forall u\in \LieGOtimesaKtq,~f\in\Ktq,$$
and hence it is indeed a $\Kt$-module morphism. To see that $\Omega$
takes values in ${L^2}$, it suffices to prove that
$[\Omega(u),f]=0$, $\forall u\in \LieGOtimesaKtq$, ~$f\in \Ktq$. In
fact,
\begin{eqnarray*}
[\Omega(u),f]&=&
[d_*u-[\Lambda,u],f]\\
&=&d_*[u,f]-[u,d_*f]-[[\Lambda,f],u]-[\Lambda,[u,f]]\\
&=&0.
\end{eqnarray*}
Moreover, since  both $d_*$ and $[\Lambda,\cdot] $ are derivations,
so is $\Omega$. In other words, it is a 1-cocycle. We claim that
$(\Lambda,\Omega)$ is a compatible pair. In fact, we have the
identity
$$
d_*^2(u)=[\frac{1}{2}[\Lambda,\Lambda]+\Omega(\Lambda)
,u]+\Omega^2(u)\;,\;\quad\forall u\in \Ktq\otimes\LieG . $$ Note
that this equation already implies that
$$
d_*^2(f)=[\frac{1}{2}[\Lambda,\Lambda]+\Omega(\Lambda)
,f],\quad\forall f\in\Ktq.
$$
Therefore the compatibility of the pair is equivalent to $d_*^2=0$.

We show that if two compatible pairs $(\Lambda,\Omega)$ and
$(\Lambda',\Omega')$ correspond to the same Lie Rinehart bialgebra
$(\LieGOtimesaKtq,d_*)$, then they are equivalent. In fact,  from
the assumption, we have
$$d_*f=[\Lambda,f]=[\Lambda',f],\;\quad\forall f\in
\Ktq.$$ Hence $\Lambda'-\Lambda\in {L^2}$. We set $\Lambda'-\Lambda
= \nu$ and it follows that $\Omega' = \Omega-[\nu,\cdot]$.

Conversely, given a compatible pair $(\Lambda,\Omega)$, then
$(\LieGOtimesaKtq,d_*)$ is clearly a Lie Rinehart bialgebra, where
the operator $d_*:~\Ktq\otimes \wedge^k\LieG \lon \Ktq\otimes
\wedge^{k+1}\LieG $ is defined by formula
(\ref{EqtRelationDStarAndPair}).
\end{proof}

\vskip 0.2cm

By (3) of Proposition \ref{Pro:Existencegamma} and the above proof,
we have the following conclusion.
\begin{thm}\label{Thm:compatiblepairdstar2}
  Let $\theta$
be an action of $\LieG$ on $\Kt$   of Type 3. Then there is a
one-to-one correspondence between Lie Rinehart bialgebras
$(\LieGOtimesaKt ,d_*)$   and equivalence classes of compatible
pairs $(\Lambda,\Omega)$   such that $d_*=[\Lambda,\cdot]+\Omega$.
Here $\Lambda\in \Kt\otimes \wedge^2\LieG $ and $\Omega: ~~\LieG\lon
L\wedge_{\Kt} L$, $L$ being the kernel of $\theta:
\LieGOtimesaKt\lon \Kt$.
\end{thm}

Most of our paper concentrate on non-trivial actions. We now examine
the case where  the action is zero.  Let $\cawu$ be the ring $\Kt$
or $\Ktq$. For $Y\in \LieGOtimesa$ we have an induced derivation
$d_Y$:~$\cawu\otimes \wedge^{k} \LieG \lon \cawu\otimes \wedge^{k+1}
\LieG  $,
$$
d_Y(fW)\defbe f' Y\wedge W,\quad\forall f\in \cawu, W\in
\wedge^k\LieG.
$$
Note that $d_Y$ is not able to be written as $[\Lambda,\cdot]$, for
some $\Lambda\in \cawu\otimes \wedge^2\LieG$.

We omit the proof of the following proposition.
\begin{prop}For the trivial crossed product $\LieGOtimesa$, any Lie
Rinehart bialgebra  $(\LieGOtimesa, d_*)$ is uniquely determined by
a $1$-cocycle $\Omega$ and an element $Y\in \cawu\otimes Z(\LieG)$
such that
$$
d_*=\Omega+d_Y.
$$
Moreover, $\Omega$ and $Y$ are subject to the following
conditions:
$$ \Omega(Y)=0,\quad \Omega^2+d_Y\circ\Omega=0,\mbox{ as
a map }\LieG\lon \cawu\otimes\wedge^2 \LieG.
$$
\end{prop}

Especially when $\LieG$ is semisimple, we will show that $d_*$ is
always a coboundary  in the next section (Theorem
\ref{Thm:trivialactionbialgebra}).

\section{Lie Rinehart Bialgebras for Semisimple $\Kt$-Crossed Products}
\label{Sec:SemisimpeactionLRBaglebras} In this section, we  study
the special properties of a crossed product $\LieGOtimesaKt$ and a
Lie Rinehart bialgebra $(\LieGOtimesaKt,d_*)$, where $\LieG$ is
semisimple.

\begin{thm}\label{Thm:trivialactionbialgebra}
Let $\LieG$ be a semisimple Lie algebra and $\cawu$ be an arbitrary
algebra. For the trivial crossed product  $\LieGOtimesa$ (i.e.,
$\theta=0$), the Lie Rinehart bialgebras $(\LieGOtimesa,d_*)$ are
one-to-one in correspondence with $\Lambda\in \cawu\otimes
\wedge^2\LieG $ satisfying $[X,[\Lambda,\Lambda]]=0$, $\forall
X\in\LieG$, such
that $d_*=[\Lambda,\cdot]$.  
\end{thm}

We will need the famous Whitehead's lemma which claims that for  any
nontrivial, finite dimensional $\LieG$-module $V$, the cohomology
groups $\huaH^1(\LieG,V)$ and $\huaH^0(\LieG,V)$ are both zero
\cite{MR0143793}.

\begin{proof}[Proof of Theorem \ref{Thm:trivialactionbialgebra}] Since
$\LieGOtimesa$ is a freely generated $\cawu$-module, one is easy to
see that $(\LieGOtimesa,d_*)$ becomes a Lie Rinehart bialgebra if
and only if all the following three conditions hold.
\begin{eqnarray}\label{Eqt:dStarOnXY}
d_*[X,Y]&=&[d_*X,Y]+[X,d_*Y],\\\label{Eqt:dStarOnaY}
d_*[a,Y]&=&[a,d_*Y]+[d_*a,Y],\\\label{Eqt:dStarOnab}
[d_*a,b]&=&-[d_*b,a],
\end{eqnarray}
$\forall X,Y\in \LieG, a,b\in \cawu$.

It is quite evident that $d_*(\LieG)$ is contained in a subspace
$\frkD=\sum_{i=1}^m a^i (\wedge^2\LieG)$, $a^i\in \cawu$, $m\in
\mathbb{N}$, which is clearly a $\LieG$-module. It follows from
relation (\ref{Eqt:dStarOnXY}) that   $d_*|_{\LieG}$ is a
$1$-cocycle, and by the   Whitehead's lemma, there exists
$\Lambda\in \frkD$ such that $d_*|_{\LieG}=[\Lambda,\cdot]$.

Notice that for the trivial crossed product  $\LieGOtimesa$,
condition (\ref{Eqt:dStarOnaY}) becomes $[d_*a,Y]=0$, $\forall Y\in
\LieG$. Then again by the Whitehead's lemma, we know that $d_*a=0$,
$\forall a\in \cawu$. Thus, $d_*(x)=[\Lambda,x]$ holds for all $x\in
\LieGOtimesa$. Clearly, the condition $[X,[\Lambda,\Lambda]]=0$ is
equivalent to $d_*^2(X)=0$.

The uniqueness of $\Lambda$ is guarantied by the Whitehead's lemma.
\end{proof}

\vskip 0.2cm

In what follows, we will study a nontrivial  crossed product
$\Kt\otimes\LieG $, where $\frkg$ is a semisimple Lie algebra
possessing a nontrivial action $\theta$ on $\Kt$ of Type 3. We will
classify all Lie Rinehart bialgebras $(\LieGOtimesaKt,d_*)$.

 By Theorem \ref{Thm:sl2KOtimesKtAllEqual},
 we know that $\frkg$
must be  of the form: $\frkg = \splinear(2,K) \oplus \frkl$ where
$\frkl = \Ker\theta$ is  an arbitrary semisimple Lie algebra. By
means of the  Killing form $\textbf{(}\cdot, \cdot \textbf{)}$ of
$\frkg$, one can identify $\frkg^*$ with $\frkg$  and define the
Cartan 3-form $\Omega$ by
$$\Omega(X,Y,Z)= ([X, Y], Z), ~~~\forall X, Y, Z \in \frkg,$$
which is a Casimir element $\Omega\in \wedge^3\LieG$ (i.e.,
$[\Omega,X]=0$, $\forall X\in \LieG$).  In particular, we denote the
Cartan 3-form of $\splinear(2,K)$ by $\Omega_{sl(2)}$. Under the
base $ \EEM ,\HH,\EEP $ of $\splinear(2,K)$ given in
(\ref{Eqt:HEPEMofsl2}), the values of the Killing forms are
determined by
\begin{equation}\label{Killingofsl2}
(\HH,\HH)=2, \ (\EEP ,\EEM )=(\EEM ,\EEP )=-4. \end{equation}
Therefore, we have $\Omega_{sl(2)}=4\HH\wedge\EEP \wedge\EEM $.

The Killing form is naturally extended to be a product of
$\LieGOtimesaKt$, taking values in $\Kt$. For each $f\in \Kt$, we
denote $(d_\theta f)^{\#}\in \LieGOtimesaKt$   the corresponding
element for $d_\theta f \in \LieGXingOtimesaKt$, i.e.,
$$
((d_\theta f)^{\#},X)=\theta(X)f,\quad\forall X\in \LieG.
$$
We introduce a differential operator  from $\Kt\otimes \wedge^k\LieG
$ to $ \Kt\otimes \wedge^{k+1}\LieG$ as follows,
\begin{equation}\label{Diff}
\Diff( f X_1\wedge\cdots\wedge X_k ) =( d_\theta f)^{\#}\wedge
X_1\wedge\cdots\wedge X_k, ~~\forall f\in \Kt,~~~~ \forall X_1,
\cdots, X_k\in \LieG.
\end{equation}
The operator $\Diff$ is totally determined by $\Diff t$ since $\Diff
f = f^{' }\Diff t,~~\forall f \in \Kt$ and $\Diff X = 0,~~ \forall X
\in \LieG$.
\begin{lem}\label{lem:diff1-2}
 \begin{equation}\label{Eqt:DiffSquaret}
\Diff^2 t=\frac{1}{32}[\Omega_{sl(2)},t].
\end{equation}
\end{lem}
\begin{proof} For the  standard $\theta$ given in Example
\ref{Ex:sl2KstandardMorphism}, we have
$$(\Diff t, \mathbf{E}_i)=\theta(\mathbf{E}_i)=t^i,\ \ i=0,1,2.
$$
Thus  the relations in (\ref{Killingofsl2}) implies
\begin{equation}\label{Eqt:Difft1}
\Diff t=(d_\theta t)^{\#}=\frac{1}{4}(2t \HH-t^2\EEM -\EEP ),
\end{equation}
and
\begin{eqnarray}\nonumber
\Diff^2 t &=& \frac{1}{4}(2\Diff t \wedge\HH-2t\Diff t\wedge\EEM ),
\\\nonumber
&=&\frac{1}{8}(t\EEP \wedge\EEM +E_1\wedge\EEP - t^2E_1\wedge \EEM
)\\\nonumber &=&\frac{1}{8}[\HH\wedge\EEP \wedge\EEM ,t].
\end{eqnarray}
The latter one is exactly $\frac{1}{32}[\Omega_{sl(2)},t]$. By
Theorem \ref{Thm:sl2KOtimesKtAllEqual}, this relation must hold for
any nontrivial $\theta$. \end{proof}

\begin{defn}\label{Def:epsilonDYBE}
With notations above, for a constant $\varepsilon$ and  an element
$\Lambda  \in \Kt\otimes\wedge^2\LieG $, the following equation is
called the $\varepsilon$-dynamical Yang-Baxter equation
($\varepsilon$-DYBE):
\begin{equation}\label{Eqt:cYBE}
\half[\Lambda,\Lambda]+\varepsilon\Diff \Lambda +
\frac{\varepsilon^2}{32}\Omega_{sl(2)} = \omega \in (\wedge
^3\frkl)^{\frkl},
\end{equation}
where $\omega$ is an arbitrary  Casimir element in $\wedge ^3\frkl$.
A solution to this equation is called an $\varepsilon$-dynamical
$r$-matrix.
\end{defn}
We remark that this notion is a special one of the notion of
dynamical $r$-matrices coupled with Poisson manifolds introduced in
\cite{MR1930079}, which is a natural generalization of the classical
dynamical $r$-matrices of Felder \cite{MR1404026}.

The main theorem in this section is as follows:
\begin{thm}\label{Thm:LieGOtimesKtBialgebraStrucuture}
 For any Lie Rinehart algebra
 $\LieGOtimesaKt$, where $\frkg$ is   a semisimple Lie algebra  possessing  a
nontrivial action  on $\Kt$,
   there is a one-to-one correspondence between Lie Rinehart bialgebras  $(\LieGOtimesaKt, d_*)$
   and  $\varepsilon$-dynamical
$r$-matrices $\Lambda$ such that
$$d_*=[\Lambda,\cdot]+\varepsilon\Diff.$$
\end{thm}

We split the proof into several lemmas.

\begin{lem}\label{Thm:H1LieGWedge2LieGIsZero}
For any $K$-linear operator $D: \LieG\lon \Kt\otimes \wedge^2\LieG $
satisfying
\begin{equation}\label{Eqt:temp5}
D[X,Y]=[DX,Y]+[X,DY], \ \ \forall X,Y\in \LieG,
\end{equation}
there exists a unique  $\Lambda\in\Kt\otimes\wedge^2\LieG $ such
that $D=[\Lambda,\cdot]$.
\end{lem}
\begin{proof} Suppose that $D(X)=\sum_{i=0}^m t^iD_i(X)$, for each $X\in
\LieG$, where the operators $D_i: \LieG\lon \wedge^2\LieG$ are all
$K$-linear and $m\in \mathbb{N}$ is the highest degree appeared in
the image of $D$.

\noindent\emph{Claim 1.} $D_m(\HH)=0$. This is seen by comparing the
highest term on both sides of the relation
\begin{eqnarray*}
D(\EEP )&=&D([\HH,\EEP ])=[D(\HH),\EEP ]+[\HH,D(\EEP )]\\
&=&\sum_{i=0}^mt^i([D_i(\HH),\EEP ]+[\HH,D_i(\EEP )]+iD_i(\EEP ))
-\sum_{i=1}^mit^{i+1}D_i(\HH).
\end{eqnarray*}

\noindent\emph{Claim 2.} $D_m(\EEM )=0$. This comes from the
relation
$$-2D(\HH)=D([\EEP ,\EEM ])=[D(\EEP ),\EEM ]+[\EEP ,D(\EEM )].$$

\noindent\emph{Claim 3.} $D_m(X)=0$, $\forall X\in \kerthetazero$.
This is by $[X,\EEP ]=0$.

\noindent{\emph{Claim 4.}} $m\neq 1$. In fact, if $m=1$, we suppose
that $D_1(\EEP )=a\HH\wedge\EEP +b\HH\wedge\EEM +c\EEP \wedge\EEM $,
for some $a,b,c\in K$. Then comparing the two sides of the relation
below {\emph{Claim 1}}, one is able to get $[\HH, D_1(\EEP )]=0$,
which implies $a=b=0$. By comparing the relation below {\emph{Claim
2.}}, one gets $[\EEM , D_1(\EEP )]=0$, which implies $c=0$. Thus
$D_1(\EEP )=0$, contradicts with our assumptions that $m=1$ is the
highest degree appeared in the image of $D$.

Now, we know that $D_m(\EEP )\neq 0$. If $m\geq 2$, we define a new
operator
$$D^{(1)}\defbe D-\frac{1}{m-1}[t^{m-1}D_m(\EEP ),\cdot].$$
It obviously satisfies a 1-cocycle condition similar to
(\ref{Eqt:temp5}).  Assume that
$D^{(1)}=\sum_{i=1}^nD^{(1)}_it^i(\cdot)$, where $D^{(1)}_i:
\LieG\lon \wedge^2\LieG$ are all $K$-linear and $n$ is the highest
degree appeared in $\Img (D^{(1)})$, then clearly $n\leq m$. But it
is easily seen that
$$D^{(1)}_m(\HH)=D^{(1)}_m(\EEM )=D^{(1)}_m(\EEP )=D^{(1)}_m(\kerthetazero)=0,$$
and hence $n<m$.

In this way, the induction goes forward and it amounts to prove that
$D^{(l)}$ is a coboundary, for sufficiently large $l\in \mathbb{N}$.
It suffices to assume that $\Img (D^{(l)})\in \wedge^2\LieG$, in
which case the Whitehead's Lemma is valid and this proves that $D$
is a coboundary.

Next we show that $\Lambda$ is unique, i.e., If any $\tau\in
\Kt\otimes \wedge^2\LieG $ satisfies
$[X,\tau]=0$, $\forall X\in \LieG$,
then it must be zero.  Write $\tau=\sum_{i=0}^m t^i \tau_i$, for
some $\tau_i\in \wedge^2\LieG$ ($\tau_m\neq 0$), then $[\tau,\EEP
]=0$ becomes
$$[\tau_0,\EEP ]+\sum_{i=1}^m
t^i([\tau_i,\EEP ]-(i-1)\tau_{i-1})-m t^{m+1}\tau_m=0.$$ Thus, $m$
must be zero, $\tau\in \wedge^2\LieG$. The conclusion $\tau=0$ comes
from the fact that $\huaH^0(\LieG,\wedge^2\LieG)=0$, since $\LieG$
is semisimple.\end{proof}

\begin{rmk}This lemma suggests that   $\huaH^i
(\LieG,~\Kt\otimes \wedge^2\LieG )=0$ ($i = 1, 2$) is also true.
\end{rmk}
By Lemma \ref{Thm:H1LieGWedge2LieGIsZero}, we know that for any
1-degree derivation $d_*$ for the Gerstenhaber  algebra
$\Kt\otimes\wedge^{\bullet}\LieG $, there exists a unique
$\Lambda\in \Kt\otimes \wedge^2\LieG $ such that $d_*|_{\frkg} =
 [\Lambda,\cdot]$. The next lemma gives some further information
 on $d_*$ as follows.

\begin{lem}\label{lem:DefofDiffanddStar}
With notations above, then, for the following   operator:
$$\Diffo\defbe d_*-[\Lambda,\cdot]~:~
\Kt\otimes \wedge^k\LieG \lon \Kt\otimes \wedge^{k+1}\LieG ,$$ there
exists a constant number $\varepsilon$ such that $\Diffo =
\varepsilon \Diff$.
\end{lem}

\begin{proof} Recall the three conditions listed in the proof of Theorem
\ref{Thm:trivialactionbialgebra}. In particular, $ d_*=
[\Lambda,\cdot] + \Diffo $, which naturally subjects to
(\ref{Eqt:dStarOnXY}), is a derivation for Lie brackets if and only
if $d_*$ satisfies the other two conditions (\ref{Eqt:dStarOnaY})
and (\ref{Eqt:dStarOnab}), i.e., $[\Diffo t,t] = \theta(\Diffo t)=0
$, and $\Diffo [X,t]=[X,\Diffo t],\ \forall X\in \LieG$. Thus
$[X,\Diffo t]=0$, $\forall X\in \frkl$ and we know that $\Diffo t\in
\Kt\otimes\wedge^2 \splinear(2,K) $. Suppose that
$$\Diffo t= \alpha \HH + \beta \EEP  + \gamma \EEM ,$$
for some  $\alpha,\beta,\gamma\in \Kt$. Then one obtains
\begin{eqnarray*}
&&\Diffo t=\Diffo [\HH,t]=[\HH,\Diffo t]\\
&=&t\alpha' \HH +(t\beta'+\beta) \EEP +(t\gamma'-\gamma)\EEM .
\end{eqnarray*}
Hence  $t\alpha'=\alpha$, $beta'=0$ and $gamma'=2\gamma$. So we get
$\alpha=a t$, $\beta=b$, $\gamma=ct^2$, where $a,b,c$ are some
constants. On the other hand, we have
\begin{eqnarray*}
&&2t\Diffo t   = \Diffo t^2\\
&=&\Diffo [\EEP ,t]=[\EEP ,\Diffo t]\\
&=&(t^2\alpha'-2\gamma)\HH + (\beta'-\alpha) \EEP  + \gamma' \EEM .
\end{eqnarray*}
Hence we get
$$2t\alpha=t^2\alpha'-2\gamma,\ 2t\beta=\beta'-\alpha.
$$
These two relations restrain that $a:b:c=-2:1:1$. This proves that
there exists $\varepsilon\in K$ such that
$$\alpha=\half\varepsilon t;\ \ \beta= -\frac{1}{4}\varepsilon; \ \
\gamma=-\frac{1}{4}\varepsilon t^2.$$ Then by formula
(\ref{Eqt:Difft1}), $\Diffo t= {\varepsilon}\Diff t$. \end{proof}

\begin{lem}\label{Pro:Casimir}
For any $\Gamma\in \Kt\otimes \wedge^3 \LieG $ satisfying $[\Gamma,
X]=0,\ \ \forall X\in \LieG$,  $\Gamma$ must be of the form
$\Gamma=k\Omega_{sl(2)}+\omega$, where $k$ is a constant  and
$\omega$ is a Casimir element in $\wedge^3\frkl$.
\end{lem}
\begin{proof} Using the same method as in the proof of Theorem
\ref{Thm:H1LieGWedge2LieGIsZero}, one easily gets $\Gamma \in
\wedge^3\LieG$. So we write $$\Gamma=\Gamma^{3,0}+
\Gamma^{2,1}+\Gamma^{1,2}+\Gamma^{0,3},$$ where $\Gamma^{ij}\in
\wedge^i\splinear(2,K)\wedge (\wedge^j\frkl)$. Clearly,
$\Gamma^{0,3}$ is a Casimir element in $\wedge^3\frkl$, and so is
$\Gamma^{3,0}$. If we write $\Gamma^{1,2}=\HH\wedge A+\EEP \wedge
B+\EEM \wedge C$, where $A,B,C\in\wedge^2\frkl$, then
$[\Gamma,\frkl]=0$ implies $[A,\frkl]=0$. Since $\frkl$ is
semisimple, $A$ must be zero. Similarly, $B=C=0$, and
$\Gamma^{1,2}=0$. For the same reasons, $\Gamma^{2,1}=0$.
\end{proof}

\begin{proof}[Proof of Theorem
\ref{Thm:LieGOtimesKtBialgebraStrucuture}] By Lemma
\ref{lem:DefofDiffanddStar}, any 1-degree derivation $d_*$ for the
Gerstenhaber  algebra $\Kt\otimes\wedge^{\bullet}\LieG $ has a
unique decomposition
$$d_*=[\Lambda,\cdot]+\varepsilon\Diff,$$
where one does not need any compatible conditions between $\Lambda$
and $\Diff$. It is easy to check that $d^2_* = 0$ if and only if
\begin{equation}\label{Eqt:thetaPairFirst}
 [\half[\Lambda,\Lambda]+\varepsilon\Diff \Lambda,t]+\varepsilon^2\Diff^2
t =0
\end{equation}
and
\begin{equation}\label{Eqt:thetaPairSecond}
[\half[\Lambda,\Lambda]+\varepsilon\Diff \Lambda,X]= 0,\ \ \forall
X\in \LieG.
\end{equation}

Now, due to   (\ref{Eqt:thetaPairSecond}) and Lemma
\ref{Pro:Casimir}, we have
$$\half[\Lambda,\Lambda]+\varepsilon\Diff
\Lambda=k\Omega_{sl(2)}+\omega.
$$ Moreover, by
(\ref{Eqt:thetaPairFirst}), (\ref{Eqt:DiffSquaret}), we obtain
$k=-\frac{\varepsilon^2}{32}$. That is exactly (\ref{Eqt:cYBE}).
\end{proof}

\begin{cor}\label{Cor:FormulasofDual}
Identifying $\Kt\otimes\LieG^* $  with $\LieGOtimesaKt$ via the
Killing form,  for the second (dual) Lie Rinehart algebra structure
on $\LieGOtimesaKt$, the Lie bracket and the anchor, are given by
the following formulas,
\begin{equation}\label{Eqt:XingBracketFormula}
[x,y]_*=[x,y]_{\Lambda}+\varepsilon (\theta(x).y-\theta(y).x)\,, \
\forall x,y\in \LieGOtimesaKt.
\end{equation}
and
\begin{equation}\label{Eqt:ThetaXingFormula}
\theta_*=\theta\circ(\Lambda^\sharp+\varepsilon I).
\end{equation}
Moreover, under these two structures, $\LieGOtimesaKt$ is a crossed
product  if and only if $\Lambda\in\wedge^2\LieG$.
\end{cor}
\begin{proof} It is some straightforward calculations to verify formulas
(\ref{Eqt:XingBracketFormula}) and (\ref{Eqt:ThetaXingFormula}). In
particular, for $X,Y\in \LieG$, by relation
(\ref{Eqn:phipsiLambdabracket2}), we have
$$([X,Y]_*,Z)=([X,Y]_\Lambda,Z)=([Z,\Lambda],X\wedge Y),
\ \forall Z\in \LieG.$$ Thus, $[X,Y]_*\in\LieG$, holds for all $X$,
$Y$ if and only if $[\LieG,\Lambda]\in\wedge^2\LieG$, which simply
suggests $\Lambda\in \wedge^2\LieG$. Only when this happens,
$\LieGOtimesaKt$ endowed with the dual bracket and anchors, becomes
a crossed product.\end{proof}

\begin{prop}\label{Pro:eteqwt;eral;jkg}
There exists some $\tau\in \Kt\otimes\wedge^2 \LieG $ such that
\begin{itemize}

\item[1)] $\Diff t=[\tau,t]$ and $\tau$ is unique up to  an
element of $\wedge^2_{\Kt} L$, where $L$ is the kernel of $\theta:
\LieGOtimesaKt\lon \Kt$.

\item[2)] The operator defined by $\Omega\defbe
\Diff-[\tau,\cdot]$, $E\lon \wedge_{\Kt}^2 L$ is a $1$-cocycle with
respect to the adjoint representation.

\item[3)] One  can take such $\tau\in \Kt\otimes \wedge^2
\splinear(2,K) $ which is also an $\varepsilon$-dynamical $r$-matrix
for
 $\varepsilon = -1$.
\end{itemize}
\end{prop}

\begin{proof} We first prove (3). Let $\theta$ be the standard action. We
can check that
$$
\tau=-\frac{1}{4}\EEP \wedge\EEM +\frac{t}{2}\HH\wedge \EEM ,
$$
satisfies $\Diff t=[\tau,t]$ (c.f. Equation (\ref{Eqt:Difft1})), and
it is a $(-1)$-dynamical r-matrix. 
This shows the existence of $\tau$ in (1).
If $\widetilde{\tau}$ is another one, then
$[\tau-\widetilde{\tau},f]=0$, $\forall f\in \Kt$ implies that
$\tau-\widetilde{\tau}\in \wedge_{\Kt}^2 L$. For the operator
$\Omega$ defined in (2), it already satisfies condition
(\ref{Eqt:Cocycle}). Then from $\Diff f=f' \Diff
t=f'[\tau,t]=[\tau,f]$, $\forall f\in \Kt$, we get
$$\Omega(fx)=\Diff f\wedge_{\Kt} x+f \Diff(x)-[\tau,f]\wedge_{\Kt} x
-f[\tau,x]=f\Omega(x),\ \ \forall x\in \LieGOtimesaKt.$$ This shows
that $\Omega$ is  a $\Kt$-linear map. \end{proof}

Now, we can determine the compatible pair declared by Theorem
\ref{Thm:compatiblepairdstar2}. In fact, the above proposition
claims that for a Lie Rinehart bialgebra $(\sltwoOtimesa, d_*)$,
$d_*=[\Lambda,\cdot]+\varepsilon\Diff$   can be written into the
form
$$
d_*=[\Lambda+\varepsilon\tau,\cdot]+\varepsilon(\Diff-[\tau,\cdot])=
[\Lambda+\varepsilon\tau,\cdot]+\varepsilon\Omega.
$$

So $(\Lambda+\varepsilon\tau,\varepsilon\Omega)$ is a compatible
pair.

It is seen that the case that $\frkg = \splinear(2,K)$ is the most
important  case,  which we shall now examine. Let $E=\Kt\otimes
\splinear(2,K)  $ be the
 Lie Rinehart algebra coming from  the standard  action
$\theta: (\HH,\EEP ,\EEM )$ $\mapsto $ $(t,t^2,1)$. Set
$$\Lambda=u
\HH\wedge\EEP + v\EEP \wedge\EEM +w \HH\wedge \EEM ,  ~~\, \, ~~~
u,v,w\in \Kt.$$ By  some straightforward calculations, one gets
$$[\Lambda,\Lambda]=
(-v^2-uw+\half t^2[u,v]+\half [v,w]+\half t[u,w])\Omega_{sl(2)}\,,$$
and
$$\Diff\Lambda=\frac{1}{16} (2tv'+w'-t^2u')\Omega_{sl(2)}\,.$$
Thus, we see that
$$ \half[\Lambda,\Lambda]+\varepsilon\Diff \Lambda   = f_{\varepsilon}\Omega_{sl(2)},$$
where  function $ f_{\varepsilon}$ is defined by
$$f_\varepsilon(u,v,w)\defbe
-\half(v^2+uw)+\frac{1}{4}(t^2[u,v]+[w,v]+t[u,w])+\frac{\varepsilon}{16}(2tv'+w'-t^2u').
$$
Consequently, we have
\begin{cor}\label{sl2KOtimesKtBi}
 Let $\Kt\otimes \splinear(2,K) $ be the
 Lie Rinehart algebra with the standard  action. Then
$$\Lambda=u
\HH\wedge\EEP + v\EEP \wedge\EEM +w \HH\wedge \EEM $$ is an
$\varepsilon$ dynamical $r$-matrix if and only if
$f_\varepsilon(u,v,w)=-\frac{1}{32}\varepsilon^2$, i.e.,
\begin{equation}\label{Eqt:CompatibleConditionForfghc}
-16(v^2+uw)+8(t^2[u,v]+[w,v]+t[u,w])+2\varepsilon(2tv'+w'-t^2u')
+\varepsilon^2=0.
\end{equation}
\end{cor}

\begin{ex}\label{Ex:check} Assume that $u=0$, $v=v_0$, $w(t)=w_0t$ where
$v_0$, $w_0$ are all constants, then
(\ref{Eqt:CompatibleConditionForfghc}) becomes
$\varepsilon^2+2w_0\varepsilon-8v_0(w_0+2v_0)=0$.
The two solutions are $\varepsilon=4v_0$ and $\varepsilon=-2w_0-4 v_0$. 
\end{ex}
\begin{ex}  Check that $u=a_0$, $v=a_0t+\frac{\varepsilon}{4}$,
$w=-\frac{\varepsilon}{2}t-a_0t^2$ is also a solution to
(\ref{Eqt:CompatibleConditionForfghc}), where $a_0$ is a constant.
\end{ex}


When $\Lambda$ belongs to $\wedge^2 \splinear(2,K)$, or $u,v,w$ are
all constants, Equation (\ref{Eqt:CompatibleConditionForfghc})
becomes $v^2+uw=\varepsilon^2/16$. So we conclude from Corollary
\ref{Cor:FormulasofDual} that
\begin{cor}
For the Lie Rinehart bialgebra $(E=\sltwoOtimesa,d_*)$, if the
induced Lie Rinehart algebra $\EXinga$ is also a crossed product,
then there exists a unique quadruple $(u,v,w,\varepsilon)\in K^4$
satisfying
\begin{equation}\label{Eqt:FourConstants}
v^2+uw={\varepsilon^2}/{16},
\end{equation} and
\begin{equation}\label{Eqt:dStar}
d_*=[u \HH\wedge\EEP + v\EEP \wedge\EEM +w \HH\wedge \EEM
,~\cdot]+\varepsilon\Diff.
\end{equation}
Conversely, any quadruple $(u,v,w,\varepsilon)\in K^4$ satisfying
{\rm(\ref{Eqt:FourConstants})} corresponds to a Lie Rinehart
bialgebra $(E,d_*)$ by relation \rm{(\ref{Eqt:dStar})} and $\EXinga$
is also a crossed product.
\end{cor}

We then consider $\LieG=\LieG_1\oplus \LieG_2$ where
$\LieG_1\cong\LieG_2\cong \splinear(2,K)$. Suppose that $\LieG_1$
acts nontrivially  on $\Kt$ and  $\Ker (\theta)=\LieG_2$. Let
(\ref{Eqt:HEPEMofsl2}) be the standard base of $\LieG_1$, and
$(\bar{\HH},\bar{\EEP },\bar{\EEM })$ be the standard base of
$\LieG_2$. Again we assume that $\theta: (\HH,\EEP ,\EEM )\mapsto
(t,t^2,1)$.

\begin{ex}
Let $\Lambda=(t+1)\bar{\HH}\wedge\bar{\EEP }+t^2\bar{\EEP
}\wedge\bar{\EEM }+ (1-t)\bar{\HH}\wedge\bar{\EEM }$ be an element
in $\Kt\otimes\wedge^2\LieG_2 $. Then, $\Lambda$ is a $0$-dynamical
$r$-matrix.
\end{ex}
Suppose that a bisection of $\Kt\otimes \wedge^2\LieG $ given by
$$\Lambda=a \HH\wedge\EEP +b\EEP \wedge\EEM +c \HH\wedge\EEM
+u\bar{\HH}\wedge\bar{\EEP }+v\bar{\EEP }\wedge\bar{\EEM }+
w\bar{\HH}\wedge\bar{\EEM },$$ where $a,b,c,u,v,w\in \Kt$. Then
\begin{eqnarray*}
&&\half[\Lambda,\Lambda]+\varepsilon\Diff \Lambda\\
&=& f_{\epsilon}(a,b,c)\Omega_{ sl(2)}-\half(v^2+uw)\bar{\Omega_{
sl(2)}}
\\&& + ((at^2+c+\half\varepsilon t)\HH+((b-\frac{1}{4}\varepsilon
)-at)\EEP -((b+\frac{1}{4}\varepsilon)t^2+ct)\EEM ) \\ &&\, \, \,
\wedge (u'\bar{\HH}\wedge\bar{\EEP }+v'\bar{\EEP }\wedge\bar{\EEM }+
w'\bar{\HH}\wedge\bar{\EEM }).
\end{eqnarray*}
Hence, $\Lambda$ is a solution to the $\varepsilon$-DYBE if and only
if $f_\varepsilon(a,b,c)=-\frac{1}{32}\varepsilon^2$, $v^2+uw$ is a
constant and
\begin{eqnarray*}
\left\{
\begin{array}{r@{\quad = \quad}l}
at^2+c+\half\varepsilon t & 0,\\
(b-\frac{1}{4}\varepsilon
)-at & 0,\\
(b+\frac{1}{4}\varepsilon)t^2+ct & 0.
\end{array} \right.
\end{eqnarray*}
There are many solutions to  these conditions. For  an example,
$a=a_0$, $b=a_0t+\frac{\varepsilon}{4}$,
$c=-\frac{\varepsilon}{2}t-a_0t^2$ ($a_0\in K$), $u=t+1$, $v=t^2$,
$w=t-1$.

\end{document}